\documentclass[10pt, oneside]{amsart}
\usepackage{geometry} 
\usepackage{amsthm, amsmath, amssymb, enumerate, enumitem, color, hyperref}
\usepackage[all]{xy}
\usepackage{graphicx}
\usepackage[parfill]{parskip}    

\DeclareMathOperator{\dist}{dist}

\DeclareMathOperator{\Q}{\mathcal{Q}}

\DeclareMathOperator{\F}{\mathcal{F}}
\DeclareMathOperator{\G}{\mathcal{G}}
\DeclareMathOperator{\id}{\textup{id}}
\DeclareMathOperator{\into}{\hookrightarrow}

\DeclareMathOperator{\Lim}{\underrightarrow{\lim}}

\DeclareMathOperator{\coker}{coker}

\newcounter{number}[section]

\newenvironment{nummer}{\refstepcounter{number}{\noindent\arabic{section}.\arabic{number}}}{}

\newcommand{\bn}{\noindent \begin{nummer} \rm}
\newcommand{\en}{\end{nummer}}

\newenvironment{ntheorem}{\noindent {\sc Theorem:} \it}{}
\newenvironment{nlemma}{\noindent {\sc Lemma:} \it}{}
\newenvironment{nprop}{\noindent {\sc Proposition:} \it}{}
\newenvironment{ndefn}{\noindent {\sc Definition:} \it}{}
\newenvironment{ncor}{\noindent {\sc Corollary:} \it}{}

\newenvironment{nproof}{\noindent {\sc Proof:}}{\mbox{}\hfill 
\rule[-.2ex]{.25em}{1.8ex}}

\title[Products of Cantor sets and odd spheres]{$\mathrm{C}^*$-algebras of minimal dynamical systems of the product of a Cantor set and an odd dimensional sphere}
\author{Karen R.~Strung}

\address{Mathematisches Institut\\
Westf\"alische Wilhelms-Universit\"at\\
Einsteinstra\ss e 62\\
48149 M\"unster\\
Germany}

\email{karen.strung@uni-muenster.de}

\date{\today}
\thanks{Supported by SFB 878, GIF Grant 1137-30.6/2011}

\begin{document}

\begin{abstract}
Let $\beta : S^{n} \to S^{n}$, for $n = 2k + 1, k \geq 1$, be one of the known examples of a nonuniquely ergodic minimal diffeomorphism of an odd dimensional sphere. For every such minimal dynamical system $(S^{n}, \beta)$ there is a Cantor minimal system $(X, \alpha)$ such that the corresponding product system $(X \times S^n, \alpha \times \beta)$ is minimal and the resulting crossed product $\mathrm{C}^*$-algebra $\mathcal{C}(X \times S^n) \rtimes_{\alpha \times \beta} \mathbb{Z}$ is tracially approximately an interval algebra (TAI). This entails classification for such \mbox{$\mathrm{C}^*$-algebras}. Moreover, the minimal Cantor system $(X, \alpha)$ is such that each tracial state on $\mathcal{C}(X \times S^n) \rtimes_{\alpha \times \beta} \mathbb{Z}$ induces the same state on the $K_0$-group and such that the embedding of $\mathcal{C}(S^n) \rtimes_{\beta} \mathbb{Z}$ into $\mathcal{C}(X \times S^n) \rtimes_{\alpha \times \beta} \mathbb{Z}$ preserves the tracial state space. This implies $\mathcal{C}(S^n) \rtimes_{\beta} \mathbb{Z}$ is TAI after tensoring with the universal UHF algebra, which in turn shows that the $\mathrm{C}^*$-algebras of these examples of minimal diffeomorphisms of odd dimensional spheres are classified by their tracial state spaces.
\end{abstract}
\maketitle

\parskip1ex
\parindent0.8em

\setcounter{section}{-1}

\section{Introduction}

The classification programme for $\mathrm{C}^*$-algebras aims to classify simple separable unital nuclear $\mathrm{C}^*$-algebras by their Elliott invariants, an invariant which can be assigned to any unital $\mathrm{C}^*$-algebra consisting of $K$-theory, the tracial state space, and a map which pairs the tracial state space with the $K_0$-group. The programme was initiated by George Elliott with his classification of the approximately finite-dimensional (AF) $\mathrm{C}^*$-algebras. He later conjectured that all simple unital nuclear $\mathrm{C}^*$-algebras should be classified by these invariants, that is, that an isomorphism of Elliott invariants should be liftable to a $^*$-isomorphism of $\mathrm{C}^*$-algebras.  

It is now known that the conjecture does not hold in full generality: various counterexamples have been constructed of nonisomorphic  simple separable unital nuclear $\mathrm{C}^*$-algebras which nevertheless are indistinguishable at the level of their Elliott invariants.  A short explanation of the failure of the conjecture can be attributed to the construction of the Jiang--Su algebra $\mathcal{Z}$ by Xinhui Jiang and Hongbing Su \cite{JiaSu:Z}. The Jiang--Su algebra is a simple separable unital nuclear $\mathrm{C}^*$-algebra with no nontrivial projections that is infinite-dimensional yet has the same Elliott invariant as $\mathbb{C}$. The conjecture then predicts that for any $\mathrm{C}^*$-algebra $A$ which falls within the scope of the programme we should have that $A$ is $\mathcal{Z}$-stable, that is,  $A \cong A \otimes \mathcal{Z}$. Counterexamples thus far produced all show failure of $\mathcal{Z}$-stability. 

The study of $\mathrm{C}^*$-algebras has long been influenced by the subject of dynamical systems.  The case where $(X, \alpha)$ is a minimal dynamical system is of particular interest to the classification programme for $\mathrm{C}^*$-algebras, as the crossed product $\mathcal{C}(X) \rtimes_{\alpha} \mathbb{Z}$ is an elegant example of a simple nuclear $\mathrm{C}^*$-algebra.  In the case of a minimal homeomorphism of a Cantor set, one has the striking classification result of Ian Putnam, Thierry Giordano and Christian Skau which shows an isomorphism of $K$-theory of the associated $\mathrm{C}^*$-algebras implies the existence of a $^*$-isomorphism of the $\mathrm{C}^*$-algebras themselves. Moreover, in this setting, isomorphisms of $K$-theory---hence also $^*$-isomorphisms of $\mathrm{C}^*$-algebras---imply topological strong orbit equivalence of the underlying dynamical systems \cite{GioPutSkau:orbit}. 

For spaces of higher dimension, classification becomes more complicated.  While many of the techniques used to tackle the examples coming from Cantor systems can be generalized to an arbitrary infinite compact metrizable space $X$, the fact that $X$ is no longer completely disconnected means that there are not necessarily many projections in the $\mathrm{C}^*$-algebra.  Nevertheless, substantial progress has been made in this more general setting, for example \cite{LinPhi:MinHom, EllEva:irrrot, LinMat:CantorT1, LinMat:CantorT2, Sun:CantorTorus}. A particularly wide-reaching result follows from the work of Andrew Toms and Wilhelm Winter in \cite{TomsWinter:minhom, TomsWinter:PNAS} (which uses a special case of a more general theorem found in \cite{StrWin:Z-stab_min_dyn}), where they give classification for the class of $\mathrm{C}^*$-algebras of minimal dynamical systems $(X, \alpha)$ of infinite finite-dimensional metric spaces under the additional assumption that projections in the $\mathrm{C}^*$-algebras separate their tracial states. 

The correspondence between tracial states on $\mathcal{C}(X) \rtimes_{\alpha} \mathbb{Z}$ and $\alpha$-invariant Borel probability measures on $X$ means that the class covered by Toms and Winter's classification includes those $\mathrm{C}^*$-algebras associated to uniquely ergodic minimal dynamical systems. However, notable examples of nonuniquely ergodic minimal dynamical systems have been constructed, and moreover, the resulting $\mathrm{C}^*$-algebras need not have projections separating tracial states. Perhaps the most famous such examples come from minimal diffeomorphisms of $n$-spheres, with $n \geq 3$ odd.  In \cite{FatHer:Diffeo}, Albert Fathi and Michael Herman give an argument that proves the existence of uniquely ergodic diffeomorphisms of $S^n$. Generalizing this result (based on the so-called ``fast approximation-conjugation'' technique of Dmitri Anosov and Anatole Katok), Alistair Windsor shows that one can construct a minimal diffeomorphism on $S^n$ with any prescribed number (finite, countable or continuum) of ergodic measures \cite{Wind:not_uniquely_ergo}.  In all these cases, Alain Connes shows that the resulting $\mathrm{C}^*$-algebras have no nontrivial projections \cite[Section 5]{Con:Thom}. Thus as soon as one has more than one ergodic measure the associated $\mathrm{C}^*$-algebras lie beyond current classification results. 

Previous classification techniques for the $\mathrm{C}^*$-algebras of minimal dynamical systems have used a large $\mathrm{C}^*$-subalgebra, originally introduced by Ian Putnam in \cite{Putnam:MinHomCantor}, which, in a sense, breaks the orbit of the homeomorphism at a single point in the underlying space. Let $(X, \alpha)$ be a minimal dynamical system and $u$ the canonical unitary implementing $\alpha$ in the $\mathrm{C}^*$-algebra $\mathcal{C}(X) \rtimes_{\alpha} \mathbb{Z}$. For any point $x \in X$ one defines a $\mathrm{C}^*$-subalgebra $\mathrm{C}^*(\mathcal{C}(X), u\mathcal{C}_0(X \setminus \{x\})) \subset \mathcal{C}(X) \rtimes_{\alpha} \mathbb{Z}$ (that is, the $\mathrm{C}^*$-subalgebra generated by $\mathcal{C}(X)$ and $u \mathcal{C}_0(X \setminus \{x\}))$. This subalgebra retains a lot of information about $\mathcal{C}(X) \rtimes_{\alpha} \mathbb{Z}$ while having a significantly more tractable structure.

In the case of the minimal Cantor systems, these subalgebras turn out to be AF algebras. More generally, N. Christopher Phillips and Qing Lin show that for arbitrary minimal dynamical systems of infinite compact metrizable spaces these subalgebras can be written as inductive limits of so-called recursive subhomogeneous (RSH) subalgebras \cite[Section 3]{LinQPhil:KthoeryMinHoms}. In the case where projections separate tracial states, then after tensoring $\mathrm{C}^*(\mathcal{C}(X), u\mathcal{C}_0(X \setminus \{x\}))$ with the universal UHF algebra $\mathcal{Q}$,  they are moreover tracially approximately finite-dimensional (TAF), as defined by Huaxin Lin \cite{Lin:TAF1}. This is then shown to pass to $\mathcal{C}(X) \rtimes_{\alpha} \mathbb{Z}$ and entails classification. Unfortunately, this technique is so far insufficient in the case of the odd dimensional spheres. For these examples, the invariant suggests that they should be tracially approximately interval algebras (TAI) after tensoring with $\Q$.  Using the results of the author and Winter in \cite{StrWin:Z-stab_min_dyn}, it would be enough to show that the large $\mathrm{C}^*$-subalgebra is TAI after tensoring with $\Q$, as this would then pass to the whole $\mathrm{C}^*$-algebra. However, at present, showing that an approximately RSH algebra is TAI after tensoring with $\Q$ requires certain restrictions which are unlikely to be satisfied in this case \cite{StrWin:slice}.

Let $(S^n, \beta)$ be a minimal dynamical system as constructed by Windsor in \cite{Wind:not_uniquely_ergo}. In this paper, instead of tensoring our $\mathrm{C}^*$-algebra with a UHF algebra, we consider the related dynamical system given by the product action with a minimal Cantor set. For a minimal system $(S^n, \beta)$, we show that we can find a minimal Cantor system $(X, \alpha)$ such that $(X \times S^n, \alpha \times \beta)$ is also minimal. Instead of considering the large subalgebra given by breaking the orbit of $\alpha \times \beta$ at a point, we break the action at a fibre $\{x \} \times S^n$. This subalgebra is TAI and we show this implies the $\mathrm{C}^*$-algebra $\mathcal{C}(X \times S^n) \rtimes_{\alpha \times \beta} \mathbb{Z}$ is also TAI.  Since the subalgebra breaks the orbit at a fibre, most of the difficult work, such a Berg's technique, takes place within the setting of the Cantor minimal system $(X, \alpha)$.   In this way, after taking the product with the Cantor system, the crossed product is easier to handle.  This entails classification of these product systems.  The structure of the arguments given are mostly based on the results of Huaxin Lin and Hiroki Matui's work on crossed products of minimal dynamical systems on $(X \times \mathbb{T})$ in \cite{LinMat:CantorT1, LinMat:CantorT2} as well as the work of H. Lin and Phillips on orbit-breaking subalgebras in \cite{LinPhi:MinHom}.

Furthermore, we are able to ensure that the embedding $\mathcal{C}(S^n) \rtimes_{\beta} \mathbb{Z} \into \mathcal{C}(X \times S^n) \rtimes_{\alpha \times \beta} \mathbb{Z}$ preserves the tracial state space and that each tracial state in $T( \mathcal{C}(X \times S^n) \rtimes_{\alpha \times \beta} \mathbb{Z})$ induces the same state on $K_0( \mathcal{C}(X \times S^n) \rtimes_{\alpha \times \beta} \mathbb{Z})$.  It  then follows from \cite[Theorem 4.2]{Win:ClassCrossProd} that $(\mathcal{C}(S^n) \rtimes_{\beta} \mathbb{Z}) \otimes \Q$ is in fact TAI. This implies classification, by tracial state spaces, for the $\mathrm{C}^*$-algebras of minimal diffeomorphisms of odd dimensional spheres constructed in \cite{Wind:not_uniquely_ergo}. That this should be true was conjectured by Phillips \cite[Conjecture 4.7]{Phi:CancelSRDirLims}.

Finally, we provide some discussion about how these results might be extended from minimal dynamical systems $(S^n, \beta)$ to more general systems $(Y, \beta)$.  Here the author would like to note that after the paper became available as a preprint, H. Lin was able to use this technique and some of the suggested generalisations to successfully cover further minimal dynamical systems \cite{Lin:Spheres}, indicating that the process of taking a product with a Cantor system has excellent potential for further applications to the classification programme. 

The author would like to thank Taylor Hines, Bhishan Jacelon, Snigdhayan Mahanta and Wilhelm Winter for helpful discussions, and the referee for pointing out an error in the initial version. Part of this work was carried out while the author was a Postdoctoral Fellow at the Fields Institute in Toronto, Canada.

\section{Breaking the orbit at a fibre}

\bn
Let $(X, \alpha)$ be a Cantor dynamical system and $\mathcal{C}(X) \rtimes_{\alpha} \mathbb{Z}$ the associated $\mathrm{C}^*$-algebra with canonical unitary $u$. A  {Kakutani--Rokhlin partition} is a family of nonempty clopen subsets of $X$, index by a finite subset $V$, denoted 
\[ \mathcal{P} = \{ X(v, k) \mid v \in V, 1 \leq k \leq h(v) \}, \]
satisfying the following: 
\begin{enumerate}
\item $\mathcal{P}$ partitions $X$, 
\item for every $v \in V$ and $k= 1, 2, \dots, h(v) - 1$ we have $\alpha(X(v,k)) = X(v, k+1)$.
\end{enumerate}
The clopen subset $R(\mathcal{P}) := \cup_{v \in V} X(v, h(v))$ is called the roof set of $\mathcal{P}$ (see \cite{HerPutSkau:Bratteli}).

If $y \in X$ then, as shown in \cite{Putnam:MinHomCantor} (see also \cite[Section VIII.6]{Dav:C*-ex} for the construction), we may find a sequence $(\mathcal{P}_k)_{k \in \mathbb{N}}$ of Kakutani--Rokhlin partitions,  $\mathcal{P}_k = \{ X(k,v,j)  \mid v \in V_k, j = 1, \dots, h_k(v) \}$ where the union generates the topology of $X$ and such that the roof sets $R(\mathcal{P}_k)$ are nested decreasing clopen sets shrinking to the singleton $\{y\}$. Such a partition is used to show that
\[ \mathrm{C}^*(\mathcal{C}(X), u\mathcal{C}_0(X \setminus \{y\})) \cong \Lim A_k\]
where the $A_k$ are AF algebras given by 
\[A_k = \oplus_{v \in V_k} M_{h_k(v)}(\mathcal{C}(X_{k,v, h_{k(v)}})).\]
\en

\bn In what follows, $n$ will always be an odd number with $n \geq 3$.

 Let $\alpha : X \to X$ and $\beta : S^n \to S^n$ be homeomorphisms and consider the product homeomorphism $\alpha \times \beta : X \times S^n \to X \times S^n$ given by 
\[ \alpha \times \beta (x,y) = (\alpha(x), \beta(y)). \]
 We will denote 
\[ A  = \mathcal{C}(X \times S^n) \rtimes_{\alpha \times \beta} \mathbb{Z}. \]
Let $u$ denote the canonical unitary implementing $\alpha \times \beta$, that is, the unitary $u$ satisfying $u f u^* = f \circ (\alpha \times \beta)^{-1}$ for every $f \in \mathcal{C}(X \times S^n)$. Note that we have canonical embeddings 
\[ \mathcal{C}(X) \rtimes_{\alpha} \mathbb{Z} \into A\]
and 
\[ \mathcal{C}(S^n) \rtimes_{\beta} \mathbb{Z} \into A. \]
\en 

\bn For a clopen subset $Y \subset X$, we denote 
\[ A_Y = \mathrm{C}^*(\mathcal{C}(X \times S^n), u\mathcal{C}_0((X \setminus Y) \times S^n)).\]  

\begin{nprop} \label{AH}  Let $y \in X$. Then $A_{\{y\}}$ is an inductive limit of AH algebras with no dimension growth.
\end{nprop}

\begin{nproof} The proof is similar to \cite[Proposition 3.3 (1)]{LinMat:CantorT1} after replacing $\mathbb{T}$ with $S^n$, based on \cite[Section 3]{Putnam:MinHomCantor} (also \cite[Section VIII.6]{Dav:C*-ex}).
For every $k\in \mathbb{N}$, let 
\[ \mathcal{P}_k = \{ X(k,v,j)  \mid v \in V_j, j = 1, \dots, h_k(v) \}, \]
be Kakutani-Rokhlin partitions corresponding to the minimal Cantor system $(X, \alpha)$ such that the roof sets 
\[ \textstyle R(\mathcal{P}_k) = \bigcup_{v \in V_j} X(k, v, h_k(v)) \]
shrink to the singleton $\{y\}$ and the union of $\mathcal{P}_k$ generate the topology on $X$. Let 
\[ A_k = \mathrm{C}^*(\mathcal{C}(X \times S^n), u \mathcal{C}_0((X \setminus R(\mathcal{P}_k) \times S^n) .\]

Since the roof sets shrink to $\{y\}$, we have 
\[ A_{\{y\}} = \Lim  A_k. \]

For $i, j \in \{1, \dots, h_k(v)\}$ define $e^{(k,v)}_{i,j} = u^{i-j}1_{X(k, v, j) \times S^n}$. It is straightforward to check that by the choice of the partition, each $e^{(k,v)}_{i,j} \in A_k$.  One checks that the $e^{(k,v)}_{i,j}$ satisfy the relations of matrix units and as in \cite[Section 3]{Putnam:MinHomCantor}, the $\mathrm{C}^*$-subalgebra
\[ \mathrm{C}^*(e^{(k,v)}_{i,j}, \{ (f, 1_{S^n}) \mid f|_{X(k, v, j)} \in \mathbb{C} \})  \cong M_{h_k(v)}\]
is finite-dimensional. Furthermore,
\[ \oplus_{v \in V_k} \oplus_{j=1}^{h_k(v)} e^{(k,v)}_{i, i} A_k e^{(k,v)}_{i, i} \cong \oplus_{v \in V_k} \oplus_{j=1}^{h_k(v)} \mathcal{C}(X(k,v,j) \times S^n),\]
hence  $A_k$ is a direct sum of homogeneous $\mathrm{C}^*$-algebras with topological dimensional most $n$ whence $A_{\{y\}} = \Lim  A_k$ is an AH algebra with no dimension growth.
 \end{nproof}
\en

\bn
\begin{nprop}
Suppose that $\alpha :X \to X$ and $\beta : S^n \to S^n$ are minimal homeomorphisms. Then the $K$-theory of $A$ is given by
\[ K_i(A) \cong \mathbb{Z} \oplus K_0(\mathcal{C}(X) \rtimes_{\alpha} \mathbb{Z})\]
for $i=0$ and $i=1$, where 
\[ K_0(\mathcal{C}(X) \rtimes_{\alpha} \mathbb{Z}) \cong \mathcal{C}(X, \mathbb{Z}) / \{f - f\circ\alpha^{-1} \mid f \in \mathcal{C}(X, \mathbb{Z})\}. \]
\end{nprop}

\begin{nproof}
By the K\"unneth Theorem for tensor products,  we have
\[ K_0(\mathcal{C}(X \times S^n)) \cong \mathcal{C}(X, \mathbb{Z}), \quad K_1(\mathcal{C}(X \times S^n)) \cong \mathcal{C}(X, \mathbb{Z}) .\]
Since $n$ is odd and $\beta: S^n \to S^n$ has no fixed points, its degree is $\text{deg}(\beta)=(-1)^{n+1} = 1$ \cite[Exercise V.1.3]{OutRuiz:MapDeg} thus $\beta$ is homotopic to the identity map by the Hopf Theorem \cite[Theorem V.2.1]{OutRuiz:MapDeg}. It follows that the homeomorphism $\alpha \times \beta$ induces the map $\alpha_i$ on the $K_i(\mathcal{C}(X \times S^n))$ given by $\alpha_i(f) = f \circ \alpha$ for $f \in \mathcal{C}(X, \mathbb{Z})$. The Pimsner--Voiculescu six term exact sequence is
\begin{displaymath}
\xymatrix{ \mathcal{C}(X , \mathbb{Z}) \ar[r]^{\id - \alpha_0} &  \mathcal{C}(X , \mathbb{Z}) \ar[r]^{\iota_0} & K_0(A)  \ar[d] \\
K_1(A) \ar[u] & \mathcal{C}(X, \mathbb{Z}) \ar[l]^{\iota_1} & \mathcal{C}(X, \mathbb{Z}) \ar[l]^{\id - \alpha_1}.
}
\end{displaymath}
Let $f \in \ker (\id - \alpha_i)$, where $i = 0$ or $1$. Then $f (x) = f \circ \alpha(x)$ for all $x \in X$. Since $\alpha$ is minimal, this implies that $f$ is constant, thus $\ker(\id - \alpha_i) \cong \mathbb{Z}$. We also have $\mathcal{C}(X, \mathbb{Z})/\{f - f\circ \alpha^{-1}  \mid f \in \mathcal{C}(X, \mathbb{Z})\}) \cong \coker(\id - \alpha_i)$. Since $\mathbb{Z} = \ker(\id - \alpha_{i+1})$ is free abelian, $K_i(A) = \coker(\id - \alpha_i) \oplus \ker(\id - \alpha_{i+1})$, that is, 
\[ K_i(A) \cong \mathbb{Z} \oplus K_0(\mathcal{C}(X) \rtimes_{\alpha} \mathbb{Z})\]
for $i=0$ and $i=1$.
\end{nproof}
\en

\section{Existence of minimal product homeomorphisms}

We now restrict to the case that $\alpha \times \beta$ is a minimal homeomorphism. Starting with a minimal homeomorphism $\beta : S^n \to S^n$ it is not necessarily the case that the product $\alpha \times \beta$ will be minimal, even if $\alpha$ is itself a minimal homeomorphism. The next result establishes that such an $\alpha$ indeed exists and moreover can be chosen to be uniquely ergodic.

The fact that an odometer system would suffice was suggested by Taylor Hines, who, at the same time as the writing of following proposition, also provided the author with a very similar proof. This result is probably known, but we couldn't find a reference so include the proof here.

\bn
Recall that an odometer is a minimal Cantor system defined as follows: Let $(m_i)_{i \in \mathbb{N}}$ be a sequence such that $2 \geq m_i$ divides $m_{i+1}$. Then there are homomorphisms $\mathbb{Z}/{m_{i+1}} \to \mathbb{Z}/{m_i}$ given by congruence modulo $m_i$ and the resulting inverse limit $\underleftarrow{\lim}\, \mathbb{Z}/m_j$ is a Cantor set. The odometer action is given by  $\alpha(x) = x+1$;  it is minimal and equicontinuous (see, for example, \cite{Down:Odometer}). 

A minimal dynamical system $(Y, \gamma)$ is called totally minimal if $\gamma^k : Y \to Y$ is minimal for every $k \in \mathbb{N}$  \cite[II.5.4]{Vries:ElemTopoDyn}. If $\beta : S^n \to S^n$ is minimal, then it is totally minimal. This follows from the fact that for any minimal dynamical system $(Y, \gamma)$ and any $k$, there is a partition of $Y$ into mutually disjoint clopen minimal $\gamma^k$-invariant subsets \cite[II.9.6.7]{Vries:ElemTopoDyn}. Since $S^n$ is connected, this implies there is only one nonempty $\beta^k$-invariant subset, namely $S^n$.

\begin{nprop}
Let $\beta : S^n \to S^n$ be a minimal homeomorphism. Then there is a uniquely ergodic minimal homeomorphism $\alpha : X \to X$ such that the homeomorphism $\alpha \times \beta : X \times S^n \to X \times S^n$ is minimal.
\end{nprop}

\begin{nproof} 
We show that $\beta$ is disjoint from any odometer, that is, for any odometer system $(X, \alpha)$ the product system $(X \times S^n, \alpha \times \beta)$ is minimal.  Let $(\underleftarrow{\lim}\, \mathbb{Z}/m_j, \alpha)$, $m_j \in \mathbb{N}$, $m_{j}$ divides $m_{j+1}$, $j \in \mathbb{N}$, be an odometer system. 

Set a metric $d$ on $X \times S^n$ and let $d_X$ and $d_{S^n}$ be the restrictions to $X$ and $S^n$, respectively.  

Suppose $(x_0, y_0) \in X \times S^n$. We will show that the orbit of $(x_0, y_0)$ is dense in $X \times S^n$. Let $\epsilon > 0$ and $(x, y) \in X \times Y$. Since $\alpha$ is equicontinuous, there is a $\delta >0$ such that $d_X(x_0', x_0) < \delta$ implies $d_X(\alpha^m(x_0'),\alpha^m(x_0)) < \epsilon/2$ for every $m \in \mathbb{N}$.  Take  $j$ sufficiently large to find elements $x_0', x' \in  \mathbb{Z}/m_j$ such that 
\[ d_X(x_0, x_0' ) < \delta \text{ and } d_X(x, x') < \epsilon/2.\]
Since $\mathbb{Z}/m_j$ is finite, there is a $k \in \mathbb{Z}$ such that $\alpha^k(x_0') = x'$. Since $\beta^{m_j}$ is minimal, there is $l \in \mathbb{N}$ such that 
\[d_{S^n}(\beta^{l m_j}(\beta^k(y_0)) , y ) < \epsilon.\]
Then we have
\begin{eqnarray*}
d((\alpha \times \beta)^{l m_j + k} (x_0, y_0) , (x, y) ) &<& d ( (\alpha^{l m_j + k}(x_0),\beta^{l m_j + k}(y_0 )),  (x, y))\\
&=& \max \{ d_X(\alpha^{l m_j + k}(x_0'), x) + \epsilon, \epsilon \}\\
&=&  \max \{ d_X(\alpha^{l m_j}(x'), x') + \epsilon, \epsilon \}\\
&=& \epsilon.
\end{eqnarray*}
\end{nproof}

\en

\bn
\begin{nprop} \label{simple}
Let $y \in X$ and let $\alpha \times \beta : X \times S^n \to X \times S^n$ be minimal. Then $A_{\{y\}}$ is simple.
\end{nprop}

\begin{nproof}
For a straightforward direct proof, one may proceed as in \cite[Proposition 2.5]{LinPhi:MinHom} after replacing each instance of $X$ with $X \times S^n$ and $U \subset X$ with $U \times S^n$. Alternatively, the result follow as in \cite[Proposition 3.3 (5)]{LinMat:CantorT1} by regarding $A$  as groupoid $\mathrm{C}^*$-algebras generated by the orbit equivalence relation, $A_{\{y\}}$ the groupoid $\mathrm{C}^*$-algebra generated by a dense subequivalence relation and applying \cite[Propostion 4.6]{Ren:groupoid}.
\end{nproof}
\en

\section{Tracial states}

\bn
\begin{nprop} Let $\alpha \times \beta : X \times S^n \to X \times S^n$ be a minimal homeomorphism where $(X, \alpha)$ is an odometer system. Then every tracial state $\tau \in T(A)$ comes from the product of the unique tracial state $\tau_0 \in T(\mathcal{C}(X) \rtimes_{\alpha} \mathbb{Z})$ and  a tracial state $\tau_1 \in T(\mathcal{C}(S^n) \rtimes_{\beta} \mathbb{Z})$. 
\end{nprop}

\begin{nproof} It is obvious that $\tau_0 \otimes \tau_1(f \circ (\alpha \times \beta)^{-1})  = \tau_0 \otimes \tau_1(f)$ for every $f \in \mathcal{C}(X \times S^n)$, so that $\tau_0 \otimes \tau_1 \in T(A)$, but we need to show that every tracial state must be of this form. Since tracial states  correspond to $(\alpha \times \beta)$-invariant probability measures, it suffices to show that every such measure on $X \times S^n$ comes from a product measure.

Let $\mu$ be an $(\alpha \times \beta)$-invariant Borel probability measure on $X \times S^n$. Define $\mu_0$ and $\mu_1$  by $\mu_0(B) = \mu(B \times S^n)$ for every Borel set $B \subset X$ and $\mu_{1}(B) = \mu(X \times B)$ for every Borel set $B \subset S^n$. Then $\mu_0$ and $\mu_1$ are Borel probability measures on $X$ and $S^n$ respectively. It is easy to seen that $\mu_0$ is $\alpha$-invariant and that $\mu_1$ is $\beta$-invariant. 

Consider the measure-preserving dynamical systems $(X, \mu_0)$, $(S^n, \mu_1)$ and $(X \times S^n, \mu)$. Since $\beta$ is totally minimal and $(X, \alpha)$ is an odometer, it follows from \cite[Lemma 12.3]{Down:Odometer} that $(X, \mu_0)$ and $(S^n, \mu_1)$ are disjoint as measurable dynamical systems, that is, $\mu =  \mu_0 \times \mu_1$ (cf. \cite{Fursten:Disjoint}).
\end{nproof}
\en

\bn
\begin{nprop}
Let $\alpha \times \beta : X \times S^n \to X \times S^n$ be a minimal homeomorphism where $(X, \alpha)$ is an odometer system.  Then $K_0(\mathcal{C}(X \times S^n)\rtimes_{\alpha \times \beta} \mathbb{Z}) = \mathbb{Z} \oplus \mathcal{C}(X, \mathbb{Z})/\{f - f\circ\alpha^{-1}\}$ has the strict ordering from the second coordinate and every tracial state $\tau \in T(\mathcal{C}(X \times S^n) \rtimes_{\alpha \times \beta} \mathbb{Z})$ induces the same state on $K_0(\mathcal{C}(X \times S^n) \rtimes_{\alpha \times \beta} \mathbb{Z})$.
\end{nprop}

\begin{nproof}
Since $\dim(X) = 0$ it follows that $H^1(X, \mathbb{Z}) = 0$ \cite[Lemma 35-3]{Nag:DimTh}. Since $H^1(S^n, \mathbb{Z}) = 0$ as well, from the K\"unneth Theorem  we get $H^1(X \times S^n, \mathbb{Z}) = 0$. In this case, to determine the range of a state $\tau_*$ induced by any tracial state $\tau \in T(A)$ it is enough to determine the range of $\tau_*$ on $K_0(\mathcal{C}(X \times S^n))$ \cite[Corollary 10.10.6]{Bla:newk-theory}. By the previous proposition, every $\tau \in T(A)$ is of the form  $\tau_X \otimes \tau_{S^n}$ where $\tau_X$ is the unique $\alpha$-invariant tracial state on $\mathcal{C}(X)$ and $\tau_{S^n}$ is a $\beta$-invariant tracial state on $\mathcal{C}(S^n)$. 

By \cite[Example 4.6]{Phi:CancelSRDirLims}, the range of any $(\tau_{S^n})_*$ is $\mathbb{Z}$ and any two $\tau_{S^n}, \tau_{S^n}'$ are homotopic via the map $t \mapsto [(1-t)\tau_{S^n} + t \tau_{S^n}']_*$. Thus the range of $\tau_X \otimes \tau_{S^n}$ is $\tau_X(K_0(\mathcal{C}(X \times \mathbb{Z}))) = \tau_X(K_0(\mathcal{C}(X))$ and since the order of $K_0(A)$ is determined by the states (for example, since $A$ is $\mathcal{Z}$-stable) it follows that the order on $K_0(A)$ is determined by the second coordinate. Finally, for any two $\tau_{S^n}, \tau_{S^n}'$, we clearly have $\tau_X \otimes \tau_{S^n}$ homotopic to $\tau_X \otimes \tau_{S^n}'$ by the above. It follows that  $(\tau_X \otimes \tau_{S^n})_* = (\tau_X \otimes \tau_{S^n}')_*$.
\end{nproof}
\en

\bn
\begin{nprop} \label{trace space}
Let $\alpha \times \beta : X \times S^n \to X \times S^n$ be a minimal product homeomorphism. Let $y \in X$ and let $u$ denote the canonical unitary in $A = \mathcal{C}(X \times S^n) \rtimes_{\alpha \times \beta} \mathbb{Z}$. Then every tracial state $\tau \in T(A_{\{y\}})$ satisfies \[ \tau(f \circ (\alpha \times \beta)^{-1}) = \tau(f)\]
for every $f \in \mathcal{C}(X \times S^n)$. 

In particular, the map $\tau \to \tau \circ \iota$ is a bijection between $T(A)$ and $T(A_{\{y\}})$.
\end{nprop}

\begin{nproof}
The proof is the same as the proof of \cite[Proposition 3.3 (4)]{LinMat:CantorT1} upon replacing $\mathbb{T}$ with $S^n$.
\end{nproof}
\en

\section{Main results}
In this section, we will show that the $\mathrm{C}^*$-algebras of the homeomorphisms of $S^n$ constructed by Windsor can be embedded into a simple unital tracially approximately interval algebra in a trace-preserving way. Let us briefly recall the definition for tracially approximately $\mathcal{S}$, where $\mathcal{S}$ is a given class of unital $\mathrm{C}^*$-algebras. 

 \bn
\begin{ndefn}\textup{(cf. \cite{Lin:amenable, EllNiu:tracial_approx})}
\label{TAS}
Let $\mathcal{S}$ denote a class of separable unital  $\mathrm{C}^{*}$-algebras. Let $A$ be a simple unital $\mathrm{C}^{*}$-algebra. Then $A$ is tracially approximately $\mathcal{S}$ (or $\mathrm{TA}\mathcal{S}$) if the following holds.

For every finite subset $\mathcal{F} \subset A$, every $\epsilon > 0$, and every nonzero positive element $c \in A$, there exist a projection $p \in A$ and a unital $\mathrm{C}^{*}$-subalgebra $B \subset pAp$ with $1_{B}=p$ and $B \in \mathcal{S}$ such that:
\begin{enumerate}
\item $\| pa - ap \| < \epsilon$ for all $a \in \mathcal{F}$,
\item $\dist(pap, B) < \epsilon$ for all $a \in \mathcal{F}$,
\item $1_{A}-p$ is Murray--von Neumann equivalent to a projection in $\overline{cAc}$.
\end{enumerate}
\end{ndefn}
\en

The $\mathrm{C}^*$-algebra $B$ is in the class I of interval algebras if it is of the form 
$$B = \bigoplus_{n=1}^N \mathcal{C}(X_n) \otimes M_{r_n}$$
for some $N \in \mathbb{N} \setminus \{0\}$, where $X_n = [0,1]$ or $X_n$ is a single point, and $r_n \in \mathbb{N} \setminus \{0\}$, $0 \leq n \leq N$.

\bn
 We now restrict to the case of homeomorphisms $\beta : S^n \to S^n$ that can be written as the limit of a sequence $(T_i)_{i \in \mathbb{N}}$ of homeomorphisms such that $T_i : S^n \to S^n$ has period $M_i$, each $M_i$ divides $M_{i+1}$, and
\[ \sup_{\substack{t \in S^n \\j = 1, \dots, M_i}}  | \beta^j(t) - T_i^j(t) | \to 0 \text{ as } i \to \infty.\]

In particular, this holds for the nonuniquely ergodic homeomorphisms constructed by Windsor \cite[Lemma 5]{Wind:not_uniquely_ergo}.

\begin{nlemma} \label{attach S^n}
There is an odometer system $(X, \alpha)$ such that the following holds:

For any $y \in X$, any $\epsilon > 0$, any $N_0 \in \mathbb{R}_+$, and any pair of finite sets $\F_X \subset \mathcal{C}(X)$, $\F_{S^n} \subset \mathcal{C}(S^n)$  there are  $M > N_0 \in \mathbb{N}$ and $y \in U \subset X$ a clopen subset and a partial isometry $w \in A_{\{y\}}  $ such that 
\begin{enumerate}
\item $\alpha^{-M}(U), \dots, \alpha^{-1}(U), U, \alpha(U), \dots, \alpha^M(U)$ are pairwise disjoint
\item $w^*w = 1_{U \times S^n}$ and $ww^* = 1_{\alpha^M(U)\times S^n}$,
\item $\| w a  - a w \| < \epsilon$ for every $a \in \{f \otimes 1_{S^n} \mid f \in \F_X\}  \cup \{1_{X} \otimes f \mid f \in \F_{S^n}\}$.
\end{enumerate}
\end{nlemma}

\begin{nproof} 
By the assumption on $\beta$, we have a sequence $(M_j)_{j \in \mathbb{N}} \subset \mathbb{N}$ with $M_j$ dividing $M_{j+1}$ and corresponding periodic homeomorphisms $T_j : S^n \to S^n$ with period $M_j$, $j \in \mathbb{N}$.  

Let $(X, \alpha)$ be the odometer system corresponding to the sequence $(M_j)_{j \in \mathbb{N}}$ (which exists since $M_j$ divides $M_{j+1}$). In this case, we may write a particular element $x \in X = \underleftarrow{\lim}\, \mathbb{Z}/M_j$ as $x = (x_j)_{j \in \mathbb{N}} \subset \Pi_{j \in \mathbb{N}} \, \mathbb{Z}/M_j$ such that $x_{j+1} = x_j \text{ mod } M_{j}$. 

Let $y  = (y_j)_{j \in \mathbb{N}} \subset X$.

Following the construction at the beginning of  \cite[Section 6]{Putnam:MinHomCantor}, we can find an approximation of $A_i \subset A_{\{y\}}$ by choosing increasingly fine Kakutani-Rokhlin partitions in $X$, as in Proposition~\ref{AH}. In this case, the roof sets $Z_i \subset X$ are given by the cylinder sets $Z_i = \{ x \in X \mid x_j = y_j, 0\leq j \leq i \}$ and the single first return time to $Z_i$ is $M_i$. The matrix units from Proposition~\ref{AH} are then given by $e_{k,l} = u^{k-l} 1_{\alpha^l(Z_i) \times S^n}$, $k, l = 1, \dots, M_i$. 

By assumption on $\beta$ there is a strictly decreasing sequence of positive numbers $(\epsilon_i)_{i \in \mathbb{N}}$ with $\epsilon_i \to 0$, such that for every $t \in S^n$ we have $d_{S^n}(\beta^{M_i}(t), t )< \epsilon_i$ for every $t \in S^n$, or correspondingly, we have that $d_{S^n}(\beta^{M_i}(\beta^{-M_i}(t) ), \beta^{M_i}(t)) < \epsilon_i$ for every $\beta^{-M_i}(t) \in S^n$. 

Let $\delta > 0$ be sufficiently small so that $d_X(x, x') < \delta$ implies $| f(x) - f(x')| < \epsilon/4$ for every $x, x' \in X$ and every $f \in \F_X$.

By the assumptions on $\beta$ as a limit of periodic homeomorphisms together with the choice of the odometer $(X, \alpha)$, by going sufficiently far out in the sequences of periods $\{M_i\}_{i \in \mathbb{N}}$, we can find  $M > N_0$ such that 
\[ d_X( \alpha^M(y), y ) < \delta/4,\]
and so that $d_{S^n}(t, \beta^{-M})$ is sufficiently small to have
\[ \| f \circ\beta^{-M} - f \| < \epsilon/2 ,\] 
 for every $f \in \F_{S^n}$.
 
 By continuity of $\alpha^M$ there is a clopen set $U \subset X$ containing $y$ of the form $\{ x \in X \mid x_j = y_j, 0\leq j \leq k \}$  for some $k \in \mathbb{N}$ such that 
\[ d_X(\alpha^M(x), x )< \delta/2 \text{ for every } x \in U. \] 
Shrinking $U$ if necessary, we may furthermore assume that $U, \alpha(U), \dots, \alpha^M(U)$ are pairwise disjoint and $d_X(x,x') < \delta/4$ for every $x, x' \in U$ and $x, x' \in \alpha^M(U)$.
 
Again by our assumptions on $\beta$, as remarked above, for sufficiently large $i$ we may arrange $d_{S^n}( t, \beta^{M_i}(t))$ is sufficently small so that  
\[ \| f\circ \beta^{-M} \circ \beta^{M_i} - f \circ \beta^{-M}\| = \sup_{t \in S^n} | f \circ \beta^{-M} \circ \beta^{M_i} (t) - f \circ \beta^{-M}(t) | < \epsilon/2 \]
for every $f \in \F_{S^n}$. By increasing $i$ if necessary, we may also arrange that the roof set of the $i^{\text{th}}$ Kakutani--Rokhlin partition $Z_i \subset U$.

Analogously to \cite[Section 6]{Putnam:MinHomCantor}, define $w_1 \in A_{\{y\}}$ by
\begin{eqnarray*}
 w_1 &=& e_{1,M_i} + \sum_{k=2}^{M_i} e_{k,k-1}\\
 &=& 1_{\alpha(U) \times S^n}  (1_{\alpha(Z_i) \times S^n}  u^{1 - M_i} +  u 1_{(X \setminus Z_i) \times S^n} ) . 
 \end{eqnarray*} 
 
Note that $w_1^* w_1 = 1_{U \times S^n}$ and $w_1 w_1^* = 1_{\alpha(U) \times S^n}$. 
 
Also, for $f \in \F_{S^n}$, we have
\begin{eqnarray*}
w_1 (1_{X} \otimes f) w_1^* &=& 1_{\alpha(Z_i)} \otimes f \circ \beta^{- 1} \circ \beta^{M_i} + 1_{\alpha(U\setminus Z_i) } \otimes f \circ \beta^{-1}.
\end{eqnarray*}

For $2 \leq m \leq M$, define $w_m \in A_{\{y\}}$ by
\[ w_m =  u 1_{\alpha^{m-1}(U) \times S^n}.\]

Put 
\[ w = w_{M}  w_{M-1}  \cdots  w_1.\] 
Then it is easy to check that $w^*w = 1_{U \times S^n}$ and $ww^* = 1_{\alpha^M(U) \times S^n}$. 

For $f \in \mathcal{C}(S^n)$ we have 
\begin{eqnarray*}
\lefteqn{w (1_{X} \otimes  f)  w^*} \\
&=& 1_{\alpha^M(Z_i)} \otimes f \circ \beta^{- M} \circ \beta^{M_i} + 1_{\alpha^M((U \setminus Z_i)} \otimes f \circ \beta^{-M},
\end{eqnarray*}
so for  $f \in \mathcal{F}_{S^n}$,
\begin{eqnarray*}
\|w (1_X \otimes f) - (1_X \otimes f) w \| &=&  \|w(1_U \otimes f) w^*  - 1_{\alpha^M(U)} \otimes f \|\\
&=& \|  f \circ \beta^{- M}  - f \| + \epsilon/2\\
& <& \epsilon.
 \end{eqnarray*}

Now let $f \in \mathcal{F}_X$. Then $| f(x) - f(x')| < \epsilon/4$ for every $x, x' \in Z_i$ by choice of $U$. Thus, since $Z_i$ is clopen in $X$, we can approximate $f$ by $\tilde{f} \in \mathcal{C}(X)$ which is constant on $Z_i$, equal to $f$ on $X \setminus Z_i$, and satisfies $\| f - \tilde{f} \| < \epsilon/4.$ Furthermore, since $d(\alpha^M(x), x )< \delta/2$ for every $ x\in U$ and both $U$ and $\alpha^M(U)$ have width less than $\delta/4$, we may assume that $\| \tilde{f}|_U - \tilde{f}|_{\alpha^M(U)} \| < \epsilon/2$. Since $\tilde{f}$ is constant on $Z_i$ we have $w_1 \tilde{f} w_1^* = \tilde{f} \circ \alpha^{-1}$ and therefore,
\begin{eqnarray*}
\| w {f} \otimes 1_{S^n} - {f} \otimes 1_{S^n} w \| &=& \| w \tilde{f} \otimes 1_{S^n} - \tilde{f} \otimes 1_{S^n} w\| + \epsilon/2\\
&=&  \| w \tilde{f}|_U \otimes 1_{S^n}w^*  - \tilde{f}|_{\alpha^M(U)} \otimes 1_{S^n}  \| + \epsilon/2\\
&<& \| \tilde{f} \circ \alpha^{-M} 1_{\alpha^M(U)} - \tilde{f}|_{\alpha^M(U)} \| + \epsilon/2 \\
&<& \| \tilde{f}|_{U} - \tilde{f}|_{\alpha^M(U)} \| + \epsilon/2 \\
&<& \epsilon.
\end{eqnarray*}
 \end{nproof}
 \en

 \bn
\begin{nlemma}  \label{proj in A_y} Let $\beta : S^n \to S^n$ be as above. There is an odometer system $(X, \alpha)$ such that the following holds:

 Let $y \in X$. Then, for any finite subset $\mathcal{F} \subset A = \mathcal{C}(X \times S^n) \rtimes_{\alpha \times \beta} \mathbb{Z}$ and any $\epsilon > 0$, there is a projection $p$ in $A_{\{y\}}  $ such that 
\begin{enumerate}
\item[\textup{(i)}] $\| pa- ap \| < \epsilon$ for all $a \in \mathcal{F}$, \label{almost com}
\item[\textup{(ii)}] $\dist (pap, p A_{\{y\}}  p) < \epsilon$ for all $a \in \mathcal{F}$, \label{almost contain}
\item[\textup{(iii)}] $\tau(1_{A  } - p) < \epsilon$ for all $\tau \in T(A  )$.
\end{enumerate}
\end{nlemma}

\begin{nproof} The proof is similar to \cite[Lemma 4.2]{LinPhi:MinHom}. Let $\epsilon> 0$. We may assume $\mathcal{F}$ is of the form
\[ \mathcal{F} = \mathcal{G}  \cup \{ u  \} \]
where $\mathcal{G} := \mathcal{G}_{X} \cup \mathcal{G}_{S^n}$ and $\mathcal{G}_{X}$, $\mathcal{G}_{S^n}$ are finite subsets of $\mathcal{C}(X)^1$ and $\mathcal{C}(S^n)^1$ respectively.

Let $M_0 \in \mathbb{N}$ such that $ \pi / (2 M_0) < \epsilon/4$.

By Lemma~\ref{attach S^n} applied to $\epsilon/4$, $M_0 + 1$, and the finite subsets $\bigcup_{m = 0}^{M_0} u^m \mathcal{G}_X u^{-m}$ and $\bigcup_{m = 0}^{M_0} u^m \mathcal{G}_{S^n} u^{-m}$, we find $M \in \mathbb{N}$ with $M > M_0 + 1$ and a clopen set $U \subset X$ containing $y$ such that
\[ \alpha^{-M_0}(U), \alpha^{-M_0 + 1}(U), \dots, U, \alpha(U), \dots, \alpha^M(U) \]
are all disjoint, and a partial isometry $w \in A_{\{y\}}  $ satisfying $w^*w = 1_{U\times S^n} := q_0$ and $ww^* = 1_{\alpha^M(U)\times S^n} := q_M$, and 
\[ \textstyle \| wf - fw \| < \epsilon/4 \text{ for all } f \in \bigcup_{k = 0}^{M_0} u^k \mathcal{G} u^{-k}. \]

Let 

\[ R > (M + M_0 + 1 )/ \min(1, \epsilon). \]

Shrinking $U$ if necessary, we may may assume that 
\[ \alpha^{-M_0}(U), \alpha^{-M_0 + 1}(U), \dots, U, \alpha(U), \dots, \alpha^R(U) \]
are pairwise disjoint. To apply Berg's technique, we only need the sets $\alpha^m(U)$ for $-M_0 \leq m \leq M$, however we require $R$ to be larger in order to satisfy property (iii) of the lemma. 

Define projections $q_m = 1_{\alpha^m(U) \times S^n} = u^m q_0 u^{-m}$. This gives us pairwise orthogonal projections in $A_{\{y\}}  $ for which conjugation by $u$ is the shift.

We now perform Berg's technique to splice along the pairs of indices $(-M_0, M - M_0), (-M_0 + 1, M - (M_0 - 1)), \dots, (0, M)$ to obtain a loop of length $M$ where conjugation by $u$ is approximately the cyclic shift.

For $t \in \mathbb{R}$, define 
$$v(t) = \cos(\pi t/2)(q_0 + q_M) + \sin(\pi t /2)(w - w^*),$$ 
which has matrix representation given by
$$ v(t) = \left( \begin{array}{ccc} \cos( \pi t/2) & -\sin(\pi t/ 2) \\ \sin( \pi t/2) & \cos(\pi t /2)  \end{array}  \right).$$
For each $t \in \mathbb{R}$, is a unitary in the corner $(q_0 + q_M) A_{\{y\}}(q_0 + q_M)$.

For $0 \leq k \leq M_0$ define 
$$w_k = u^{-k}  v(k/M_0)u^k.$$

Each $w_k \in A_{\{y \}}$ for  $0 \leq k \leq M_0$. The proof of this is the same as in the proof of \cite[Lemma 4.2]{LinPhi:MinHom}.

Estimating the matrix entries we get
\begin{eqnarray*}
\lefteqn{ \|  v((k+1)/M_0) - v(k/M_0) \|} \\
 & \leq& 2 | \cos (\pi(k+1)/2M_0 ) - \cos (\pi k/2M_0 )  | + 2  | \sin(\pi (k+1)/2 M_0) - \sin(\pi k/2 M_0 )  | \\
 & = &2 \pi/2 M_0  | \sin( \xi_1) | + 2 \pi/2 M_0 | \cos( \xi_2) |, \text{ for some } \xi_1, \xi_2 \in (k, k+1)  \\
 & \leq& 2 \pi /M_0 \\
 & <& \epsilon/2.
 \end{eqnarray*}
 From this it follows that 
$$ \| u w_{k+1} u^* - w_k \| \leq \| v((k+1)/M_0) - v(k/M_0) \| < \epsilon /2.$$

For $0 \leq m \leq M- M_0$ define $e_m = q_m$ and for $M- M_0 \leq m \leq M$ define projections 
\[ e_m = w_{M-m} q_{-(M-m)} w_{M-m}^*.\]

It is straightforward to check that the two definitions agree for $e_{M-M_0}$ and that $e_0 = e_M$.

For $1 \leq m \leq M- M_0$ conjugating $e_{m-1}$ by $u$ gives 
\[u e_{m-1}u^* =  u q_{m-1}u^* = q_m = e_m\] and for conjugation of $e_M$ we have \[ u e_{M} u^* =  u e_{0}u^* = u q_0 u^* = q_1 = e_1.\]

When $M - M_0 \leq m \leq M$ we have
\begin{eqnarray*}
\lefteqn{\|u e_{m-1} u^* - e_m \|} \\
&=&\| u w_{ {M-(m-1)}} q_{ {-(M-(m-1))}}w_{ {M-(m-1)}}^*u^* - w_{ {M-m}} q_{ {-(M-m)}} w_{ {M-m}}^*\| \\
&=& \| u w_{  {M-(m-1)}} u^* q_{ {-(M-m)}} uw_{  {M-(m-1)}}^*u^* - w_{ {M-m}} q_{ {-(M-m)}} w_{M-m}^*\| \\
& \leq& \| u w_{M-(m-1)}^*u^* - w_{M-m}^*\| +\| u w_{M-(m-1)} u^* - w_{M-m} \| \\
&=&2\| u w_{  {M-(m-1)}}^*u^* - w_{  {M-m}}^*\|\\
&<& \epsilon.
\end{eqnarray*}

Hence conjugation of the $e_m$ by $u$ is approximately a cyclic shift. 

By the definition of $e_m$, $1 \leq m \leq M$, it is clear that each $e_m \in A_{\{ y \}}$.  Set 
$$\textstyle e = \sum_{m = 1}^{M} e_m \quad \text{ and } p = 1 - e.$$
Then p is a projection in $A_{\{ y \}}$.

We now show that $p$ satisfies $(1) - (3)$ with respect to the set $\mathcal{F}$.

We have
\[
\| \textstyle p - u p u^* \| = \| \sum_{m = M - M_0 + 1}^{M} u  e_{m-1} u^* - e_m \| 
\]
The terms in the sum are pairwise orthogonal and have norm less than $\epsilon$, hence $\| p - u p u^* \|  < \epsilon$, proving $(1)$ for $u$. Also, we have $p \leq 1 - q_0 =  1 - 1_{U} \in A_{\{y\}}  $ so that $pup \in A_{\{y\}}  $, showing $(2)$ for $u $.

For $f \in \G$ we have $\| w  u^k f u^{-k}-  u^k f u^{-k}w \| < \epsilon/4$ and therefore 
\begin{eqnarray*}
\lefteqn{\| w_k f - f w_k\|}\\
&=& \|u^{k} v(k/M_0) u^{-k} f - f u^{k}   v(k/M_0)  u^{-k} \|\\
&=& \|  \sin(\pi k/(2M_0) )(w - w^*)u^{-k}  f u^k - u^{-k} f u^k  \sin(\pi k/(2M_0) )(w - w^*)\|\\
&=& 2 \| w u^{-k}  f u^k -u^{-k}  f u^k w \| \\
&<& \epsilon/2.
\end{eqnarray*}
Thus 
\begin{eqnarray*}
\lefteqn{\| p f - f p \| }\\
&=& \textstyle \| \sum_{m= M-M_0+1}^M (w_{M-m} q_{-(M-m)} w^*_{M-m} f - f w_{M-m} q_{-(M-m)} w^*_{M-m})\| \\
&<& \max_{mM} \|w_{M-m} q_{-(M-m)} w^*_{M-m} f - w_{M-m} q_{-(M-m)} f w^*_{M-m}\|  + \epsilon/2 \\
&<& \epsilon.
\end{eqnarray*}
Since $f \in \mathcal{C}(X \times S^n)$ we have that $pf p \in p A_{\{y\}}   p$. This shows $(1)$ and $(2)$ for $\G$.

For $(3)$, we have 
\[ \textstyle \tau(1 - p ) = \sum_{m=1}^M \tau(e_m) = \sum_{m=1}^{M} \tau(q_0)  < M/R  < \epsilon. \]
for every $\tau \in T(A  )$.
\end{nproof}
\en

\bn
\begin{nlemma} \label{sub TAI = TAI}  Let $A$ be a simple unital $\mathrm{C}^*$-algebra with strict comparison. Suppose that for every finite subset $\mathcal{F} \subset A$, every $\epsilon > 0$, and every nonzero positive $c \in A$, there exists a projection $p \in A$ and a simple unital $\mathrm{C}^{*}$-subalgebra $B \subset pAp$ which  is TAI, satisfies $1_{B}=p$ and
\begin{enumerate}
\item[\textup{(i)}] $\| pa - ap \| <  \epsilon$ for all $a \in \mathcal{F}$, 
\item[\textup{(ii)}] $\dist(pap, B) < \epsilon$ for all $a \in \mathcal{F}$,
\item[\textup{(iii)}] $1_{A} - p$ is Murray--von Neumann equivalent to a projection in $\overline{c Ac}$.
\end{enumerate}
Then $A$ is TAI.
\end{nlemma}

\begin{nproof} By \cite[Theorem 3.2]{Lin:simpleNuclearTR1} $A$ has property (SP) and by assumption, strict comparison. After noting this, the proof is essentially the same as that of \cite[Lemma 4.4]{LinPhi:MinHom} and \cite[Lemma 4.5]{StrWin:Z-stab_min_dyn}, replacing the $\mathrm{C}^*$-subalgebra of tracial rank zero (respecively TA$\mathcal{S}$) with the TAI $\mathrm{C}^*$-subalgebra $B$, and replacing the finite-dimensional (respectively in the class $\mathcal{S}$) $\mathrm{C}^*$-subalgebra with a $\mathrm{C}^*$-subalgebra from the class I.
\end{nproof}
\en

\bn
\begin{ntheorem} \label{product tai}
There is an odometer system $(X, \alpha)$ such that $\mathcal{C}(X \times S^n) \rtimes_{\alpha \times \beta} \mathbb{Z}$ is TAI.
\end{ntheorem}

\begin{nproof}
Set $A := \mathcal{C}(X \times S^n) \rtimes_{\alpha \times \beta} \mathbb{Z}$. Let $\epsilon> 0$ and nonzero $c \in A^1_+$ be given. Again we may assume $\mathcal{F}$ is of the form
\[ \mathcal{F} =  \mathcal{G} \cup \{ u\} \]
where $\mathcal{G} := \mathcal{G}_{X} \cup \mathcal{G}_{S^n}$ and $\mathcal{G}_{X}$, $\mathcal{G}_{S^n}$ are finite subsets of $\mathcal{C}(X)^1$ and $\mathcal{C}(S^n)^1$ respectively.

Since $A$ is $\mathcal{Z}$-stable \cite[Theorem B]{TomsWinter:PNAS}, $A$ has strict comparison  \cite[Theorem 3.1]{Win:ssa-Z-stable}. Thus we need only verify the conditions of Lemma~\ref{sub TAI = TAI}.

Let $y \in X$ and use Lemma~\ref{proj in A_y} with respect to $\F$ and $\epsilon_0 = \min \{\epsilon, \min_{\tau \in T(A)} \tau(c) \}$ to find a projection $p \in A_{\{y\}}$. Let $B = p A_{\{y\}} p$. Since $pA_{\{y\}}p$ is a simple AH algebra with no dimension growth it is TAI \cite{Lin:simpleNuclearTR1}. We have satisfied (i) and (ii) of Lemma~\ref{sub TAI = TAI} by the choice of $p$. Since $\tau(1 - p) < \epsilon_0$ for every $\tau \in T(A)$, condition (iii) follows from strict comparison.
\end{nproof}
\en

We are now able to use the above to apply \cite[Theorem 4.2]{Win:ClassCrossProd}  which entails classification of the minimal dynamical systems of odd dimensional spheres constructed by Windsor. 

\bn
\begin{ncor}
Suppose $\beta_1, \beta_2 : S^n \to S^n$ are minimal homeomorphisms as in \ref{attach S^n}. Then $A := \mathcal{C}(S^n) \rtimes_{\beta_1} \mathbb{Z} \cong \mathcal{C}(S^n) \rtimes_{\beta_2} \mathbb{Z}=: B$ if and only if $T(A) \cong T(B)$.

In particular, this holds if $\beta_1$ and $\beta_2$ come from the constructions given in \textup{\cite{Wind:not_uniquely_ergo}}.
\end{ncor}

\begin{nproof}
By the above, there are minimal Cantor systems $(X, \alpha_i)$ such that $\mathcal{C}(S^n) \rtimes_{\beta_i} \mathbb{Z}$ embeds in a trace-preserving way into $\mathcal{C}(X \times S^n) \rtimes_{\alpha_i \times \beta_i} \mathbb{Z}$,  $i \in \{1,2\}$. Moreover, $(\tau_i)_* = (\tau_i')_*$ for every tracial state $\tau_i \in T(\mathcal{C}(X \times S^n) \rtimes_{\alpha_i \times \beta_i} \mathbb{Z})$, $i \in \{1,2\}$.

By Theorem~\ref{product tai}, $\mathcal{C}(X \times S^n) \rtimes_{\alpha_i \times \beta_i} \mathbb{Z}$ is TAI, thus by \cite[Theorem 4.2]{Win:ClassCrossProd} $\mathcal{C}(S^n) \rtimes_{\beta_i} \mathbb{Z}) \otimes \Q$ is TAI, where $\Q$ denotes the universal UHF algebra. Since $A$ and $B$ satisfy the UCT, classification up to $\mathcal{Z}$-stability by Elliott invariants follows from \cite[Corollary 11.9]{Lin:asu-class}. Since $\dim(S^n) < \infty$, both $A$ and $B$ are $\mathcal{Z}$-stable \cite[Theorem B]{TomsWinter:PNAS} (or \cite[Theorem 0.2]{TomsWinter:minhom}). It follows from \cite[Example 4.6]{Phi:CancelSRDirLims} that $A$ and $B$ have isomorphic $K$-theory, thus the Elliott invariant collapses to the tracial state space.
\end{nproof}
\en

\section{Outlook} 

A close look at the method of proof suggests that one can generalize to other minimal systems of the form $(X \times Y, \alpha \times \beta)$ with $X$ a Cantor set. The particular requirements of $\beta : S^n \to S^n$ are that it is totally minimal and that $\beta$ is a limit of periodic homeomorphisms in the sense of \ref{attach S^n}.

Total minimality holds at least for any connected space $Y$. In fact, one does not need total minimality, but rather only requires that $(Y, \beta)$ is disjoint from some odometer system $(X, \alpha)$; in this case $\alpha \times \beta$ will be minimal and $T(\mathcal{C}(X \times S^n) \rtimes_{\alpha \times \beta} \mathbb{Z}) \cong T(\mathcal{C}(S^n) \rtimes_{ \beta} \mathbb{Z})$. This will hold in other situations, for example, if $(Y, \beta)$ is weakly mixing \cite[Theorem II.3]{Fursten:Disjoint}.

The restriction that $\beta$ is a limit of periodic homeomorphisms in the sense of \ref{attach S^n} will be the case if $\beta$, for example, is constructed by the ``fast approximation-conjugation'' technique as initiated in \cite{AnoKat:Fast} (using \cite[3.1]{AnoKat:Fast} with \cite[Lemma 3.7]{FayKat:EllipDyn}, provided the periods of the periodic homeomorphisms increase quickly enough). 

Otherwise, observe that the structure for $\beta$ is only required in two places in Lemma~\ref{attach S^n}. First, it is used to find the unitary $w$ required for Berg's technique.   One should be able to get around this by showing approximate unitary equivalence of the maps $\phi_1, \phi_2 : \mathcal{C}(S^n) \to q_0 A_{\{y\}} \otimes \Q q_0$ given by $\phi_1 (f) = v^* u^M q_0 f q_0u^{-M}v$ and $\phi_2(f) = q_0fq_0$, where $v$ is the unitary which comes from Berg's technique applied only to the Cantor minimal system, that is, $v q_0 v^* = q_M$. This is the technique applied in \cite{LinMat:CantorT1, LinMat:CantorT2}. Since  $q_0 A_{\{y\}} \otimes \Q q_0$ is TAI, there are various results of H. Lin which may be applied (see, for example,  \cite{Lin:AHHom}). Secondly, we have the technical requirement that $\| f\circ \beta^M - f\|$ can be made small for every $f$ in the given finite subset, which at present is more difficult to remove. Nevertheless, it appears the result will hold in greater generality.

\providecommand{\bysame}{\leavevmode\hbox to3em{\hrulefill}\thinspace}
\providecommand{\MR}{\relax\ifhmode\unskip\space\fi MR }
\providecommand{\MRhref}[2]{%
  \href{http://www.ams.org/mathscinet-getitem?mr=#1}{#2}
}
\providecommand{\href}[2]{#2}


\begin{thebibliography}{10}

\bibitem{AnoKat:Fast}
{Anosov, Dmitri V. and Katok, Anatole B.}, \emph{{New examples in smooth
  ergodic theory. {E}rgodic diffeomorphisms}}, {Trudy Moskov. Mat. Ob\v s\v c.}
  \textbf{{23}} ({1970}), {3---36}.

\bibitem{Bla:newk-theory}
{Blackadar, Bruce}, \emph{{$K$-Theory for Operator Algebras}}, {Second} ed.,
  {Mathematical Sciences Research Institute Publications}, vol.~5, {Cambridge
  University Press}, 1998.

\bibitem{Con:Thom}
{Connes, Alain}, \emph{{An analogue of the {T}hom isomorphism for crossed
  products of a {$C^{\ast} $}-algebra by an action of {${\bf R}$}}}, {Adv. in
  Math.} \textbf{{39}} ({1981}), no.~{1}, {31--55}.

\bibitem{Dav:C*-ex}
{Davidson, Kenneth R.}, \emph{{$C^*$-algebras by Example}}, Fields Institute
  Monographs, Amer. Math. Soc., Providence, R.I., 1996.

\bibitem{Vries:ElemTopoDyn}
{de Vries, J.}, \emph{{Elements of topological dynamics}}, {Mathematics and its
  Applications}, vol. {257}, {Kluwer Academic Publishers Group}, {Dordrecht},
  {1993}.

\bibitem{Down:Odometer}
{Downarowicz, Tomasz}, \emph{{Survey of odometers and {T}oeplitz flows}},
  {Algebraic and topological dynamics}, {Contemp. Math.}, vol. {385}, {Amer.
  Math. Soc.}, {2005}, pp.~{7--37}.

\bibitem{EllEva:irrrot}
{Elliott, George A. and Evans, David E.}, \emph{{The structure of the
  irrational rotation $C^*$-algebras}}, Ann. of Math. \textbf{138} (1993),
  no.~2, 477--501.

\bibitem{EllNiu:tracial_approx}
{Elliott, George A. and Niu, Zhuang}, \emph{{On tracial approximation}}, {J.
  Funct. Anal.} \textbf{254} (2008), no.~2, 396--440.

\bibitem{FatHer:Diffeo}
{Fathi, Albert and Herman, Michael}, \emph{{Existence de diff{\'e}omorphismes
  minimaux}}, {Ast{\'e}risque} \textbf{{49}} ({1977}), {37--59}.

\bibitem{FayKat:EllipDyn}
{Fayad, Bassam and Katok, Anatole}, \emph{{Constructions in elliptic
  dynamics}}, {Ergodic Theory Dynam. Systems} \textbf{{24}} ({2004}), no.~{5},
  {1477--1520}.

\bibitem{Fursten:Disjoint}
{Furstenberg, Harry}, \emph{{Math. Systems Theory}}, {Disjointness in ergodic
  theory, minimal sets, and a problem in {D}iophantine approximation}
  \textbf{{1}} ({1967}), {1--49}.

\bibitem{GioPutSkau:orbit}
{Giordano, Thierry and Putnam, Ian F. and Skau, Christian F.},
  \emph{{Topological orbit equivalence and $C^*$-crossed products}}, {J. Reine
  Angew. Math.} \textbf{{469}} ({1995}), {51--111}.

\bibitem{HerPutSkau:Bratteli}
{Herman, Richard H. and Putnam, Ian F. and Skau, Christian F.}, \emph{{Ordered
  Bratteli diagrams, dimension groups and topological dynamics}}, {Internat. J.
  Math.} \textbf{{3}} ({1992}), no.~{6}, {827--864}.

\bibitem{JiaSu:Z}
{Jiang, Xinhui and Su, Hongbing}, \emph{{On a simple unital projectionless
  {$C^*$}-algebra}}, {Amer. J. Math.} \textbf{{121}} ({1999}), no.~2,
  {359--413}.

\bibitem{Lin:amenable}
{Lin, Huaxin}, \emph{{An introduction to the classification of amenable
  {$C^*$}-algebras}}, {World Scientific Publishing Co. Inc.}, {2001}.

\bibitem{Lin:TAF1}
\bysame, \emph{{Tracially AF $C^{*}$-algebras}}, Trans. Amer. Math. Soc.
  \textbf{353} (2001), 693--722.

\bibitem{Lin:simpleNuclearTR1}
\bysame, \emph{Simple nuclear {$C^*$}-algebras of tracial topological rank
  one}, J. Funct. Anal. \textbf{251} (2007), no.~2, 601--679.

\bibitem{Lin:asu-class}
\bysame, \emph{{Asymptotic unitary equivalence and classification of simple
  amenable {$C^*$}-algebras}}, {Invent. Math.} \textbf{{183}} ({2011}),
  no.~{2}, {385--450}.

\bibitem{Lin:AHHom}
\bysame, \emph{{Homomorphisms from AH-algebras}}, arXiv preprint
  math.OA/arXiv:1102.4631v4, {2011}.

\bibitem{Lin:Spheres}
\bysame, \emph{{Minimal dynamical systems on connected odd dimensional
  spaces}}, {arXiv preprint math.OA/1404.7034}, {2014}.

\bibitem{LinMat:CantorT1}
{Lin, Huaxin and Matui, Hiroki}, \emph{{Minimal dynamical systems on the
  product of the {C}antor set and the circle}}, {Comm. Math. Phys.}
  \textbf{{257}} ({2005}), no.~{2}, {425--471}.

\bibitem{LinMat:CantorT2}
\bysame, \emph{{Minimal dynamical systems on the product of the {C}antor set
  and the circle. {II}}}, {Selecta Math. (N.S.)} \textbf{{12}} ({2006}),
  no.~{2}, {199--239}.

\bibitem{LinPhi:MinHom}
{Lin, Huaxin and Phillips, N. Christopher}, \emph{{Crossed products by minimal
  homeomorphisms}}, {J. Reine Angew. Math.} \textbf{{641}} ({2010}), {95--122}.

\bibitem{LinQPhil:KthoeryMinHoms}
{Lin, Qing and Phillips, N. Christopher}, \emph{{Ordered {$K$}-theory for
  {$C^\ast$}-algebras of minimal homeomorphisms}}, {Operator algebras and
  operator theory ({S}hanghai, 1997)}, {Contemp. Math.}, vol. {228}, {Amer.
  Math. Soc.}, 1998, pp.~{289--314}.

\bibitem{Nag:DimTh}
{Nagami, Kei{\^o}}, \emph{{Dimension theory}}, {With an appendix by Yukihiro
  Kodama. Pure and Applied Mathematics, Vol. 37}, {Academic Press}, {1970}.

\bibitem{OutRuiz:MapDeg}
{Outerelo, Enrique and Ruiz, Jes{\'u}s M.}, \emph{{Mapping degree theory}},
  {Graduate Studies in Mathematics}, vol. {108}, {American Mathematical
  Society}, {Providence, RI}, {2009}.

\bibitem{Phi:CancelSRDirLims}
{Phillips, N. Christopher}, \emph{{Cancellation and stable rank for direct
  limits of recursive subhomogeneous algebras}}, {Trans. Amer. Math. Soc.}
  \textbf{{359}} ({2007}), no.~{10}, {4625--4652}.

\bibitem{Putnam:MinHomCantor}
{Putnam, Ian F.}, \emph{{The $C^*$-algebras associated with minimal
  homeomorphisms of the Cantor set}}, Pacific J. Math. \textbf{136} (1989),
  no.~2, 329--353.

\bibitem{Ren:groupoid}
J.~Renault, \emph{{A groupoid approach to $C^*$-algebras}}, Springer-Verlag,
  1980.

\bibitem{StrWin:Z-stab_min_dyn}
{Strung, Karen R. and Winter, Wilhelm}, \emph{{Minimal dynamics and
  $\mathcal{Z}$-stable classification}}, {Internat. J. Math.} \textbf{{22}}
  ({2011}), no.~{1}, {1--23}.

\bibitem{StrWin:slice}
\bysame, \emph{{UHF-slicing and classification of nuclear {$C^*$}-algebras}},
  {to appear in J. Topology Anal.}, 2013.

\bibitem{Sun:CantorTorus}
{Sun, Wei}, \emph{{Crossed product {$C^*$}-algebras of minimal dynamical
  systems on the product of the Cantor set and the torus}}, arXiv preprint
  math.OA/1102.2801, 2011.

\bibitem{TomsWinter:PNAS}
{Toms, Andrew S. and Winter, Wilhelm}, \emph{{Minimal dynamics and the
  classification of {$C^*$}-algebras}}, {Proc. Natl. Acad. Sci. USA}
  \textbf{{106}} ({2009}), no.~{40}, {16942--16943}.

\bibitem{TomsWinter:minhom}
\bysame, \emph{{Minimal {D}ynamics and {K}-{T}heoretic {R}igidity: {E}lliott's
  {C}onjecture}}, {Geom. Funct. Anal.} \textbf{{23}} (2013), no.~{1},
  {467--481}.

\bibitem{Wind:not_uniquely_ergo}
{Windsor, Alistair}, \emph{{Minimal but not uniquely ergodic diffeomorphisms}},
  Smooth ergodic theory and its applications ({S}eattle, {WA}, 1999), Proc.
  Sympos. Pure Math., vol.~69, Amer. Math. Soc., Providence, RI, 2001,
  pp.~809--824.

\bibitem{Win:ssa-Z-stable}
{Winter, Wilhelm}, \emph{{Strongly self-absorbing {$C^*$}-algebras are
  {$\mathcal{Z}$}-stable}}, {J. Noncommut. Geom.} \textbf{{5}} ({2011}),
  no.~{2}, {253--264}.

\bibitem{Win:ClassCrossProd}
\bysame, \emph{{Classifying Crossed Products}}, {arXiv preprint
  math.OA/1308.5084}, 2013.

\end{thebibliography}
\end{document}